\documentclass[12pt]{article}
\usepackage{amsthm}
\usepackage{amsfonts}
\usepackage[T1]{fontenc}
\usepackage{epsfig}

\usepackage[lflt]{floatflt}
\usepackage{epsfig}
\usepackage{graphicx}
\usepackage{amsmath}
\usepackage{amssymb}
\usepackage{hyperref}

\usepackage{color}
\usepackage{subfigure}
\usepackage[english]{babel}

\usepackage{makeidx}
\usepackage{latexsym}

\newtheorem{theorem}{Theorem}
\newtheorem{definition}{Definition}
\newtheorem{proposition}{Proposition}
\newtheorem{lemma}{Lemma}

\newtheorem{remark}{Remark}

\newcommand{\bbN}{{\mathbb N}}

\newcommand{\bbR}{{\mathbb R}}

\newcommand{\al}{\alpha}

\newcommand{\p}{\partial}
\newcommand{\be}{\beta}
\newcommand{\G}{\Gamma}

\begin{document}
\begin{center}\emph{}
\LARGE
\textbf{Two different fractional Stefan problems which are convergent to the same classical Stefan problem}
\end{center}

                   \begin{center}
                  {\sc Sabrina D. Roscani and Domingo A. Tarzia}\\
 CONICET - Depto. Matem\'atica,
FCE, Univ. Austral,\\
 Paraguay 1950, S2000FZF Rosario, Argentina \\
(sabrinaroscani@gmail.com, dtarzia@austral.edu.ar)
                   \vspace{0.2cm}

       \end{center}
      
\small

\noindent \textbf{Abstract: }

  Two fractional Stefan problems are considered by using Riemann-Liouville and Caputo derivatives  of order $\al\in (0,1)$ such that  in the limit case ($\al=1$) both problems coincide with the same  classical Stefan problem.  For the one and the other problem, explicit solutions in terms of the Wright functions are presented. We prove that these solutions are different even though they converge, when $\al \nearrow 1$, to the same classical solution. This result also shows  that  some limits are not commutative when fractional derivatives are used.



\noindent \textbf{Note:} This paper is now published (in revised form) in Math Meth Appl Sci (2018), p.p. 1-9, DOI:10.1002/mma-5196, and is available online at \url{wileyonlinelibrary.com/journal/mma}, so always cite it with the journal's coordinates.

\section{Introduction}\label{Sec:Intro} 
\noindent In this paper two fractional Stefan problems are considered. These kind of problems are  free boundary problems where the governed equation is a  fractional diffusion equation in the temporal variable $t$. \\
A one--phase classical Stefan problem  for a semi--infinite material with initial and boundary conditions can be formulated as
\begin{equation}{\label{St-Clasico-1}}
\begin{array}{llll}
     (i)  &   \frac{\p}{\p t}u(x,t)=\lambda\, \frac{\p^2}{\p x^2} u (x,t), &   0<x<s(t), \,  0<t<T,  \, \,\\
     (ii) &   u(x,0)=f(x), & 0\leq x\leq b=s(0),\\ 
       (iii)  &  u(0,t)=g(t),   &  0<t\leq T,  \\

         (iv) & u(s(t),t)=0, & 0<t\leq T, \\
        
                 (v) & \frac{d}{dt}s(t)=-k\frac{\p}{\p x} u (s(t),t), & 0<t\leq T, 
                                             \end{array}
                                             \end{equation}
where $\lambda$ is the diffusivity  and $k$ is the conductivity of the material. This kind of problems have been widely studied (see \cite{Alexiades, Cannon, Tarzia}).\\

\noindent The fractional Caputo derivative (\cite{Ca:1967}) in the $t$ variable is defined by 
\begin{equation}\label{DefiCap} ^{C}_{0}D^{\al}_t  u (x,t)=\, _0I^{1-\al}_t u_t(x,t)= \frac{1}{\Gamma(1-\al)} \displaystyle\int^{t}_{0}\frac{\frac{\p}{\p t} u(x,\tau)}{(t-\tau)^{\al}}d\tau, \end{equation}

\noindent where $_0I^{\be}_t f(x,t)= \frac{1}{\Gamma(\be)} \displaystyle\int^{t}_{0}\frac{ f(x,\tau)}{(t-\tau)^{1-\be}}d\tau$ is the fractional Riemann--Liouville integral defined for every $\be >0$, and $\G$ is the Gamma function. \\

\noindent If we replace  in problem (\ref{St-Clasico-1}) the time derivative by the Caputo derivative  (\ref{DefiCap}), then the following fractional one--phase Stefan problem is obtained:

\begin{equation}{\label{St-Caputo}}
\begin{array}{llll}
     (i)  &  _0^CD_t^{\al} u(x,t)=\lambda^2 \frac{\p}{\p x^2}u(x,t), &   0<x<s(t), \,  0<t<T, \,  \\
     (ii) &  u(x,0)=f(x), & 0\leq x\leq b=s(0),\\ 
       (iii) &  u(0,t)=g(t),   &  0<t\leq T,  \\
         (iv) & u(s(t),t)=0, & 0<t\leq T,\\
  (v) & _0^CD^{\al}s(t)=-u_x(s(t),t), & 0<t\leq T.                                              \end{array}
  \end{equation}

 \noindent Some works foccused in problems like (\ref{St-Caputo}) are \cite{At:2012, Bl:2014, JiMi:2009, LiSun:2015, RoSa:2013, RoTa:2014, Voller:2014}. \\

\noindent Let us aboard now the physical approach. The classical mathematical model for heat flux is through the Fourier law, which says that the heat flux is proportional to the temperature gradient
 \begin{equation}\label{Fourier Law}
q^l(x,t)=-k \frac{\p}{\p x} u (x,t).
\end{equation} 
But, in the last 40 years, many generalizations of the Fourier law has been proposed  \cite{Cha:1986, HeIg:1999, IgOs:2010,  JoPr:1989, PK:1976} giving rise to the emergence of new models. In particular,  Gurtin and Pipkin \cite{GuPi:1968} propose the following law for the heat conduction, characterized by the non-locality given by
$$   q=-k\int_0^t K(t-\tau)\nabla u(\tau)d\tau,
$$
and different theories  can be developed from the consideration of different kernels of convolution. For example, in  \cite{Povstenko, Voller:2013}  a non local flow given by 
\begin{equation}\label{grad fr}  
q=-k\,^{RL}_{0}D^{1-\al}_t \frac{\p}{\p x} u (x,t)
\end{equation}
\noindent is considered, where the  fractional derivative is the Riemann--Liouville derivative respect on time of order $1-\al$ ($\al $ $\in (0,1)$)  defined by 
$$ ^{RL}_{0}D^{1-\al}_t \frac{\p}{\p x} u (x,t)= \frac{1}{\Gamma(\al)} \frac{\p}{\p t}\displaystyle\int^{t}_{0}\frac{\frac{\p}{\p x} u(x,\tau)}{(t-\tau)^{1-\al}}d\tau, \qquad \al \in (0,1).$$
\noindent Note that the non local flux coincide with the Fourier flux for $\al=1$, because $^{RL}_{0}D^{0}_t=Id$.\\

\noindent  So, we consider this non local flux.  If  (\ref{grad fr}) is replaced in the heat balance equation,  then  a fractional diffusion equation for  the fractional Riemann-Liouville derivative is  obtained:  
\begin{equation}\label{FDE}   \frac{\p}{\p t}u(x,t)=\lambda\, \frac{\p}{\p x}\left(^{RL}_{0}D^{1-\al}_t \frac{\p}{\p x} u (x,t)\right).
\end{equation}
\noindent Recalling  that $^{RL}_{0}D^{1-\al}_t$ is the left inverse operator of $_{0}I^{1-\al}_t $,  we can  apply $^{RL}_{0}D^{1-\al}_t$ to both sides of equation $(\ref{St-Caputo}-i)$  obtaining, under certain hypothesis, the  fractional diffusion equation (\ref{FDE}).\\
Fractional diffusion equations for Caputo derivatives, like $ (\ref{St-Caputo}-i)$, are linked to the modeling of diffusive processes in heterogeneous media,  such called sub or super diffusive processes (see \cite{BF:2005, Gusev, MK:2000,  Yuste:2010}).\\

\noindent Now, let us focus in the Stefan condition.  The classical Stefan condition derived in a one-phase Stefan problem is given by 
 \begin{equation}\label{cond clasica de stefan}
\frac{d}{dt}s(t)=\left.q^l(x,t)\right|_{(s(t)^-,t)}, \quad 0<t\leq T.
\end{equation}
where $q^l$ is the local flux given by (\ref{Fourier Law}). So,  replacing the non local flux (\ref{grad fr}) in (\ref{cond clasica de stefan}) we obtain the  following ``fractional Stefan condition'':
$$
\frac{d}{dt}s(t)=-\left.^{RL}_{0}D^{1-\al}_t \frac{\p}{\p x} u (x,t)\right|_{(s(t),t)}, \quad 0<t\leq T.
$$
\noindent Therefore, the second fractional Stefan problem that we can consider is given by:
\begin{equation}{\label{St}}
\begin{array}{llll}
     (i)  &   \frac{\p}{\p t}u(x,t)=\lambda\, \frac{\p}{\p x}\left(^{RL}_{0}D^{1-\al}_t \frac{\p}{\p x} u (x,t)\right), &   0<x<s(t), \,  0<t<T,  \, \,\\
     (ii) &   u(x,0)=f(x), & 0\leq x\leq b=s(0),\\ 
       (iii)  &  u(0,t)=g(t),   &  0<t\leq T,  \\

         (iv) & u(s(t),t)=0, & 0<t\leq T, \\
        
                 (v) & \frac{d}{dt}s(t)=-\left.^{RL}_{0}D^{1-\al}_t \frac{\p}{\p x} u (x,t)\right|_{(s(t),t)}, & 0<t\leq T. 
                                             \end{array}
                                             \end{equation}

\noindent  The last formulation is not usually considered because of the singularity of the Riemann--Liouville derivative, and also because the Caputo derivative is a better choice for posing fractional initial-boundary problems for  fractional parabolic operators.\\

\noindent We have seen that equations $(\ref{St}-i)$ and $(\ref{St-Caputo}-i)$ are closely linked. But, what happen with the fractional Stefan conditions $(\ref{St}-v)$ and $(\ref{St-Caputo}-v)$? \\
For example, if we apply $^{RL}_{0}D^{1-\al}_t$ to both sides of the Stefan condition $(\ref{St-Caputo}-v)$ we get
$$
\frac{d}{dt}s(t)=- ^{RL}_{0}D^{1-\al}_t \frac{\p}{\p x} u (s(t),t)
$$
which is not exactly condition $(\ref{St}-v)$, unless $\al=1$. In fact, the right side of $(\ref{St}-v)$ is
\begin{equation}\label{interch-limit}
\begin{split} 
-\left.\,^{RL}_{0}D^{1-\al}_t \frac{\p}{\p x} u (x,t)\right|_{(s(t),t)} & =-\displaystyle\lim_{x\rightarrow s(t)}\,^{RL}_{0}D^{1-\al}_t \frac{\p}{\p x} u (x,t)\\
& = -\displaystyle\lim_{x\rightarrow s(t)}\frac{\p}{\p t}\frac{1}{\G(\al)}\int_0^t(t-\tau)^{\al-1}\frac{\p}{\p x} u (x,\tau)d\tau.
\end{split}
\end{equation}

\noindent The aim of this paper is to show  explicit solutions to problems $(\ref{St-Caputo})$ and $(\ref{St})$ respectively and prove that they are different, which clearly implies that the ``fractional Stefan conditions'' $(\ref{St}-v)$ and $(\ref{St-Caputo}-v)$ are different and that for fractional derivatives some limits like (\ref{interch-limit}) are not commutative.

%$$\text{?`} \displaystyle\lim_{x\rightarrow s(t)}\frac{\p}{\p t}\frac{1}{\G(\al)}\int_0^t(t-\tau)^{\al-1}\frac{\p}{\p x} u (x,\tau)d\tau  = \frac{\p}{\p t}\frac{1}{\G(\al)}\int_0^t(t-\tau)^{\al-1}\displaystyle\lim_{x\rightarrow s(\tau)}\frac{\p}{\p x} u (x,\tau)d\tau ?  $$

\section{Previous Results}

\begin{definition} For every $x\in \bbR$ ,  \textit{Wright} function is defined by
\begin{equation}\label{W} W(x;\rho;\be)=\sum^{\infty}_{k=0}\frac{x^{k}}{k!\G(\rho k+\be)} ,\quad  \rho>-1 \text{ and }  \be\in \bbR.\end{equation}
An important particular case of of a Wright function is the \textit{Mainardi} function defined by 
$$ M_\rho (x)= W(-x,-\rho,1-\rho)=\sum^{\infty}_{n=0}\frac{(-x)^n}{n! \G\left( -\rho n+ 1-\rho \right)}, \quad \,0<\rho<1.  $$
\end{definition}

\begin{proposition}\label{Props W}  Let $\rho \in (0,1)$ be. Then the next assertions follows:
\begin{enumerate}
\item  Let $ \be\in \bbR$ be. For every $x\in \bbR$ we have
$$ \frac{\p}{\p x} W(x,\rho,\be) = W(x,\rho,\rho+\be).$$  
\item If $\be\geq 0$, then $W\left(-x,-\rho,\beta\right) $ is a positive and strictly decreasing  function in $\bbR^+$.
\item Let $\al>0$ and $ \be$ $\in \bbR$ be. For every $x>0$ and $c>0$,      
\begin{equation}\label{eq pskhu}
_0I^\al_x\,x^{\beta-1}W(-cx^{-\rho},-\rho, \beta)=x^{\beta+\al-1}W(-cx^{-\rho},-\rho, \beta+\al) .
\end{equation}
\end{enumerate}
\end{proposition}

\proof See \cite{Wr2:1940} for 1. Item 2 follows from  \cite[Theorem 8]{St:1970} and the chain rule.  Item 3 is a particular case of  \cite[Corollary 5]{Pskhu-Libro}.
\endproof

\begin{lemma}\label{gamma} For every $n\in \bbN$, it holds that {\rm \cite{BoTa:2017}}
\begin{enumerate}
\item $(2n)!=2^nn!(2n-1)!!$
\item $\G\left( n+\frac{1}{2} \right)=\frac{(2n-1)!!}{2^n}\sqrt{ \pi}.  $
\end{enumerate}
where the definition $(2n-1)!!=(2n-1)(2n-3)\cdots 5\cdot 3\cdot 1$ is used for compactness expression.
\end{lemma}

\begin{proposition}\label{conv M y W cuando al tiende a 1}Let $x\in \bbR^+_0$ be. Then  the following limits hold:\\
\begin{equation}\label{limite-M}\lim\limits_{\al\nearrow 1}M_{\al/2}\left(2x\right)=\lim\limits_{\al\nearrow 1}W\left(-2x,-\frac{\al}{2},1-\frac{\al}{2}\right)= M_{1/2}(2x)=\frac{e^{-x^2}}{\sqrt{\pi}},
\end{equation}
\begin{equation}\label{limite-W}
\lim\limits_{\al\nearrow 1}W\left(-2x,-\frac{\al}{2},\frac{\al}{2}\right)=\frac{e^{-x^2}}{\sqrt{\pi}},
\end{equation}
\begin{equation}\label{limite-erf}\lim\limits_{\al\nearrow 1}\left[1-W\left(-2x,-\frac{\al}{2}, 1\right)\right]={\rm erf}(x),
\end{equation} 
and
\begin{equation}\label{limite-erfc}\lim\limits_{\al\nearrow 1}\left[W\left(-2x,-\frac{\al}{2}, 1\right)\right]={\rm erfc}(x),
\end{equation}
where $erf(\cdot)$ is the error function defined by ${\rm erf}(x)=\frac{2}{\sqrt{\pi}}\displaystyle\int_0^xe^{-z^2}dz$ and ${\rm erfc}(\cdot)$ is the complementary error function defined by ${\rm erfc}(x)=1-{\rm erf}(x)$.\\
Moreover, the convergence is uniform over compact sets.

\end{proposition}  
\proof See \cite{RoSa:2013} for (\ref{limite-M}) and (\ref{limite-erf}).\\
 Now, for proving (\ref{limite-W}) let $\al$ be such that $0<\al< 1$. From (\ref{W}), 
\begin{equation}\label{W serie} W\left(-2x;-\frac{\al}{2};\frac{\al}{2}\right)=\displaystyle\sum^{\infty}_{k=0}\frac{(-2x)^{k}}{k!\G\left(-\frac{\al}{2} k+\frac{\al}{2}\right)}.
\end{equation}
\noindent  Let us limit the series by a convergent series which not depend on $\al$, so we can apply the Weierstrass M-test and interchange the series and the limit. \\
Recall that for all $x\, \in  \bbR,$ \cite{Erdelyi-V1}
\begin{equation}\label{Gamma-1}
\frac{1}{\G(x)\G(1-x)}=\frac{\sin (\pi x)}{\pi} ,
\end{equation} 
 and for every $k \in \bbN$
\begin{equation}\label{Gamma-2}
\G(k+1)=k\G(k).
\end{equation}
Then, 
\begin{equation}\label{Gamma-3}
\begin{split}
\left|\frac{1}{k! \G\left(-\frac{\al}{2} k+\frac{\al}{2}\right)}\right|=&\left|\frac{1}{\G(k+1) \G\left(1-\frac{\al}{2} k+\frac{\al}{2}-1\right)}\right|\\
=&\left|\frac{\G\left(\frac{\al}{2} k-\frac{\al}{2}+1\right) \sin(\pi(\frac{\al}{2} k-\frac{\al}{2}+1))}{\pi\G(k+1) }\right|\\
 \leq & \left|\frac{\G\left(\frac{\al}{2} (k-1)+1\right) }{\pi \G(k+1) }\right|.
\end{split}
\end{equation}
\noindent Now, let $x^*>0$ be the abscissa of the minimum of the Gamma function and let $k_0 $ such that $\frac{\al}{2}(k_0-1)+1>x^*$. Applying that the Gamma function is an increasing function in $(x^*,+\infty)$  it yields 
\begin{equation}\label{Gamma-4}
\left|\frac{\G\left(\frac{\al}{2} (k-1)+1\right) }{\G(k+1) }\right|\leq  \frac{\G\left(\frac{k}{2}+\frac{1}{2}\right) }{\G(k+1) }\quad \text{for all} \, k\geq k_0.
\end{equation}

\noindent Let us separate in even and odds terms. If $k=2n$, $n\in \bbN$, then applying Lemma \ref{gamma} it results that 
\begin{equation}\label{Gamma-4-1}
 \frac{\G\left(\frac{k}{2}+\frac{1}{2}\right) }{\G(k+1) }= \frac{\G\left(n+\frac{1}{2}\right) }{\G(2n+1) }= \frac{(2n-1)!!\sqrt{\pi}}{2^{n}(2n)!}< \frac{1}{2^n}=\frac{1}{\sqrt{2}^k}
\end{equation}
%\frac{\frac{(2n-1)!!}{2^n}\sqrt{ \pi}}{2^nn!(2n-1)!!}=

If $k=2n+1$, $n\in \bbN$, from Lemma \ref{gamma} we have
\begin{equation}\label{Gamma-4-2}
\begin{split}
 \frac{\G\left(\frac{k}{2}+\frac{1}{2}\right) }{\G(k+1) } &= \frac{\G\left(n+1\right) }{\G(2n+2) }=
 \frac{n!}{(2n+1)!}= \frac{n!}{(2n+1)2^{n}n!(2n-1)!!}=\\
 & =\frac{1}{(2n+1)2^{n}(2n-1)!!}\leq \frac{1}{2^{n+1}}= \frac{1}{2^{\frac{k+1}{2}}}<\frac{1}{\sqrt{2}^k}.
 \end{split}
\end{equation}
From (\ref{Gamma-4-1}) and (\ref{Gamma-4-2}) we can state that 
\begin{equation}\label{Gamma-5}
 \frac{\G\left(\frac{k}{2}+\frac{1}{2}\right) }{\G(k+1) }\leq \frac{1}{ \sqrt{2}^k}\quad \text{for all} \, k\geq k_0.
\end{equation}
\noindent From (\ref{Gamma-3}), (\ref{Gamma-4}) and (\ref{Gamma-5}) it results that the series (\ref{W serie}) is bounded by a convergent series that not depend on $\al$. Taking the limit when $\al \nearrow 1$, using  (\ref{Gamma-1}) and Lemma \ref{gamma}, the limit $(\ref{limite-W})$ holds: 

$$ \hspace{-8.5cm}\lim_{\al \nearrow 1} W\left(-2x;-\frac{\al}{2};\frac{\al}{2}\right)= $$
\begin{equation*}\begin{split}
=&\displaystyle\sum_{k=0}^{\infty} \lim\limits_{\al \nearrow 1}\frac{x^{2k}}{(2k)!\G(-\frac{\al}{2}2k+1-\frac{\al}{2})}+
\displaystyle\sum_{k=0}^{\infty}\lim\limits_{\al \nearrow 1}\frac{-x^{2k+1}}{(2k+1)!\G(1-\al(k+1))}\\
=&\displaystyle\sum_{k=0}^{\infty} \frac{x^{2k}}{(2k)!\G(-k+\frac{1}{2})}=\displaystyle\sum_{k=0}^{\infty} \frac{x^{2k}\G(k+\frac{1}{2})\sin(\pi((-k+1/2)))}{\pi(2k)!}\\
=& \frac{1}{\sqrt{\pi}}\displaystyle\sum_{k=0}^{\infty} \frac{(-x^2)^{k}}{4^kk!}= \frac{1}{\sqrt{ \pi}}e^{-\frac{x^2}{4}}.\end{split}
\end{equation*}
\endproof

%\newpage
\begin{remark} Proposition 2 shows that two different Wright functions \linebreak 
$ \G\left(1-\frac{\al}{2}\right)M_{\al/2}\left(2x\right)$ and $\G\left(\frac{\al}{2}\right)W\left(-2x,-\frac{\al}{2},\frac{\al}{2}\right) $ are convergent to the Gaussian function $G(x)=e^{-x^2}$.  A graphic for a particular value $\al=\frac{3}{4}$ is given below and 
the key  of this article is to prove that these functions does not intersect for any positive real value.

\begin{figure}[h]
\includegraphics[width=7cm]{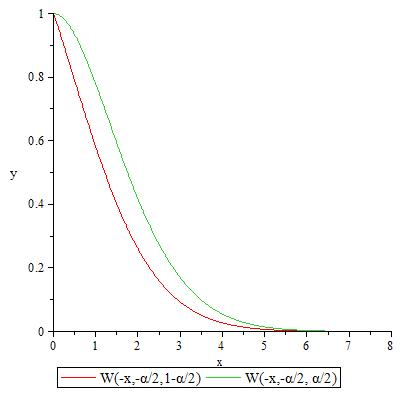}
%\caption{Figure 1}
\end{figure}
\end{remark}

\vspace{2cm}
\begin{proposition}\label{Propo monot Wrights} Let  $x>0$ be and let $0< \rho\leq \mu < \delta$. Then
$$
\G(\delta)W(-x, -\rho, \delta)<\G(\mu) W(-x,-\rho, \mu).
$$
\end{proposition}

\proof
Consider $\al=\delta-\mu$ and $\beta=\mu$ in (\ref{eq pskhu}). Then 
\begin{equation}\label{1}
_0I^{\delta-\mu}_y\,y^{\mu-1}W(-cy^{-\rho},-\rho, \mu)= y^{\delta-1}W(-cy^{-\rho},-\rho, \delta)
\end{equation}
Making the substitution $y=x^{-1/\rho}$, using (\ref{1}) and Proposition \ref{Props W}, it yields that 
\begin{equation*}
\begin{split} 
W(-x, -\rho, \delta) & =W(-y^{-\rho}, -\rho, \delta)=y^{-\delta+1}\,_0I^{\delta-\mu}_y\,y^{\mu-1}W(-y^{-\rho},-\rho, \mu)= \\
& =  y^{-\delta+1}\frac{1}{\G({\delta-\mu})} \int\limits_0^y t^{\mu-1}W(-t^{-\rho},-\rho, \mu)(y-t)^{\delta-\mu-1}dt\\
&<  y^{-\delta +1}\frac{1}{\G({\delta-\mu})} W(-y^{-\rho},-\rho, \mu)\int\limits_0^y t^{\mu-1}(y-t)^{\delta-\mu-1}dt\\
&= y^{-\delta +1}\frac{1}{\G({\delta-\mu})} W(-y^{-\rho},-\rho, \mu)\frac{\G(\mu-1+1)\G(\delta-\mu)}{\G(\delta)}y^{\delta-1} \\
& = \frac{\G(\mu)}{\G(\delta)}W(-y^{-\rho},-\rho, \mu)= \frac{\G(\mu)}{\G(\delta)}W(-x,-\rho, \mu).
\end{split}
\end{equation*}
\endproof

\section{Two different explicit solutions }
\noindent We consider two particular fractional Stefan problems:
\begin{equation}{\label{(P1)}}
\begin{array}{llll}
     (i)  &  _0^CD_t^{\al} u(x,t)= \frac{\p}{\p x^2}u(x,t), &   0<x<s(t), \,  0<t<T, \,  \\
     (ii) & s(0)=0,\\ 
       (iii) &  u(0,t)=1,   &  0<t\leq T,  \\
         (iv) & u(s(t),t)=0, & 0<t\leq T,\\
  (v) & _0^CD^{\al}s(t)=-u_x(s(t),t), & 0<t\leq T.                                              \end{array}
  \end{equation}
  and   
\begin{equation}{\label{(P2)}}
\begin{array}{llll}
     (i)  &   \frac{\p}{\p t}u(x,t)=\, \frac{\p}{\p x}\left(^{RL}_{0}D^{1-\al}_t \frac{\p}{\p x} u (x,t)\right), &   0<x<s(t), \,  0<t<T,  \, \,\\
     (ii) &  s(0)=0,\\ 
       (iii)  &  u(0,t)=1,   &  0<t\leq T,  \\

         (iv) & u(s(t),t)=0, & 0<t\leq T, \\
        
                 (v) & \frac{d}{dt}s(t)=-\left.^{RL}_{0}D^{1-\al}_t \frac{\p}{\p x} u (x,t)\right|_{(s(t),t)}, & 0<t\leq T. 
                                             \end{array}
                                             \end{equation}

%\begin{equation}{\label{St_1}}
%\begin{array}{llll}
%     (i)  &   \frac{\p}{\p t}u(x,t)=\, \frac{\p}{\p x}\left(^{RL}_{0}D^{1-\al}_t \frac{\p}{\p x} u (x,t)\right), &   0<x<s(t), \,  0<t<T,  \,\\
%      (ii)  &  u(0,t)=1,   &  0<t\leq T,  \\
%         (iii) & u(s(t),t)=0, & 0<t\leq T, \, s(0)=0 \\
       
%       (iv) & \frac{d}{dt}s(t)=-\left.^{RL}_{0}D^{1-\al}_t \frac{\p}{\p x} u (x,t)\right|_{(s(t),t)}, & 0<t\leq T, 
%                                             \end{array}
%                                             \end{equation}
%and
%\begin{equation}{\label{St_1-Caputo}}
%\begin{array}{llll}
 %    (i)  &  _0^CD_t^{\al} u(x,t)= \frac{\p}{\p x^2}u(x,t), &   0<x<s(t), \,  0<t<T, \, \\
  %          (ii) &  u(0,t)=1,   &  0<t\leq T,  \\
  %       (iii) & u(s(t),t)=0, & 0<t\leq T, \, s(0)=0 \\
  %(v) & _0^CD^{\al}s(t)=-u_x(s(t),t), & 0<t\leq T.                                              \end{array}
 % \end{equation}
\noindent It was proved in \cite{RoSa:2013} that the pair $\{w_\al,r_\al\}$ is a solution to problem \rm{(\ref{(P1)})} where 
\begin{equation}{\label{u-C}}\begin{split}
 w_\al(x,t)=&\, 1-\frac{1}{1-W\left(-2\eta_\al,-\frac{\al}{2},1\right)}\left[1-W\left(-\frac{x}{t^{\al/2}},-\frac{\al}{2},1\right)\right],\\
   r_\al(t)=&\, 2\eta_\al t^{\al/2}
\end{split}
\end{equation}
and $\eta_\al$ is the unique solution to the equation 
\begin{equation}{\label{eq_eta}} 2x\left[ 1-W\left(-2x,-\frac{\al}{2},1\right)\right]=M_{\al/2}(2x)\frac{\G(1-\al/2)}{\G(1+\al/2)}, \quad x>0.
\end{equation}

\noindent By the other side, it was proved in \cite{RoTa:2017} that the pair $\{u_\al,s_\al\}$ is a solution to problem \rm{(\ref{(P2)})}  where 
\begin{equation}{\label{u-RL}}\begin{split}
 u_\al(x,t)=& \, 1-\frac{1}{1-W\left(-2\xi_\al,-\frac{\al}{2},1\right)}\left[1-W\left(-\frac{x}{t^{\al/2}},-\frac{\al}{2},1\right)\right]\\
  s_\al(t)=& \, 2\xi_\al t^{\al/2},
  \end{split}
\end{equation}

\noindent and $\xi_\al$ is the unique solution to the equation 
\begin{equation}{\label{eq_xi}}
\begin{split}
 2x\left[ 1-W\left(-2x,-\frac{\al}{2},1\right)\right]=2xW\left(-2x,-\frac{\al}{2},1\right)+ W\left(-2x,-\frac{\al}{2},1+\frac{\al}{2}\right),\\
  x>0.
  \end{split}
\end{equation}
\noindent Looking at the similarity between solutions  (\ref{u-C}) and (\ref{u-RL}), it is natural to ask whether both are the same solution, or not. 

\begin{theorem} The explicit solutions {\rm(\ref{u-C})} to problem (\ref{(P1)}), and {\rm(\ref{u-RL})} to problem (\ref{(P2)}) are different.
\end{theorem}
\proof From \cite{Wr2:1940} we know that, for every $x\in \bbR$ the next equality holds:
\begin{equation}\label{rel wrights}
xW\left(-x,-\frac{\al}{2},1\right)+ W\left(-x,-\frac{\al}{2},1+\frac{\al}{2}\right)=\frac{2}{\al}W\left(-x,-\frac{\al}{2},\frac{\al}{2}\right).
\end{equation}
\noindent Replacing equality (\ref{rel wrights}) in (\ref{eq_xi}), we can say that the parameter  $\xi_\al$ appearing in  solution (\ref{u-RL}),  is the unique solution to the equation: 
\begin{equation}{\label{eq_xi2}} 2x\left[ 1-W\left(-2x,-\frac{\al}{2},1\right)\right]=\frac{2}{\al}W\left(-2x,-\frac{\al}{2},\frac{\al}{2}\right), \quad x>0.
\end{equation}
By the other side, we know that the parameter $\eta_\al$ which is part of the solution (\ref{u-C}) to problem (\ref{(P1)}) is the unique solution to equation (\ref{eq_eta}).  

\noindent Therefore, if we suppose that solutions (\ref{u-C}) and (\ref{u-RL}) coincides, from (\ref{eq_eta}) and (\ref{eq_xi2}) we can conclude that there exists $\nu_\al>0$ such that 
$$M_{\al/2}(2\nu_\al)\frac{\G(1-\al/2)}{\G(1+\al/2)}=  \frac{2}{\al}W\left(-2\nu_\al,-\frac{\al}{2},\frac{\al}{2}\right)$$
 or equivalently,   
$$ M_{\al/2}(2\nu_\al)\G(1-\al/2)= \G(\al/2) W\left(-2\nu_\al,-\frac{\al}{2},\frac{\al}{2}\right). 
$$
But this is a contradiction from Proposition \ref{Propo monot Wrights} and  then the thesis holds.
\endproof 

\begin{theorem} If we take the limit when $\al\nearrow 1$, the solutions (\ref{u-RL}) and (\ref{u-C})  converge to the unique solution $\{u,s\}$ to the classical Stefan problem 
\begin{equation}{\label{St-Clasico}}
\begin{array}{llll}
     (i)  &  u_t(x,t)= \frac{\p}{\p x^2}u(x,t), &   0<x<s(t), \,  0<t<T,  \\
     (ii) & u(0,t)=1,   &  0<t\leq T,  \\ 
       (iii) & u(s(t),t)=0, & 0<t\leq T, s(0)=0,\\
         (iv) & s'(t)=-u_x(s(t),t), & 0<t\leq. 
								\end{array}
  \end{equation}
\end{theorem}

\proof The unique solution to problem (\ref{St-Clasico}) is given by (see e.g. \cite{Cannon, Tarzia}),
\begin{equation}\label{u-clasica}
\begin{split}
 w(x,t) &=1-\frac{1}{\mbox{erf\,}\left(\eta_1\right)}\, \mbox{erf\,}\left(\frac{x}{2 \sqrt{t}}\right),\\
  s(t) &=2\eta_1\sqrt{t},
  \end{split}
  \end{equation} where $ \eta_1$  is the unique solution to the equation   
  \begin{equation}\label{eq eta-clasica}
  \eta_1 \,\mbox{erf}\left(\eta_1\right)=  \frac{e^{-\eta_1^2}}{\sqrt{\pi}}.
  \end{equation}
  
\noindent Note that if we take $\al=1$ in equation (\ref{eq_eta}) we recover equation (\ref{eq eta-clasica}). Now, let the sequence  $\left\{ \eta_\al\right\}_\al$ be, where $\eta_\al$ is the unique positive solution to equation $(\ref{eq_eta})$. Then
$$ 2\eta_\al=M_{\al/2}(2\eta_\al)\frac{\G(1-\al/2)}{\G(1+\al/2)}+2\eta_\al W\left(-2\eta_\al,-\frac{\al}{2},1\right).
$$
\noindent If  we define the following functions for every $x \in \bbR^+$ and $0<\al<1$: 
$$ f_\al(x)=M_{\al/2}(2x)\frac{\G(1-\al/2)}{2\G(1+\al/2)}+x W\left(-2x,-\frac{\al}{2},1\right)
$$
and 
$$
f_1(x)=\frac{e^{-x^2}}{\sqrt{\pi}}+x\mbox{erfc}\left(x\right),
$$
we have that $f_\alpha(\eta_\alpha)=\eta_\alpha$, $f_1(\eta_1)=\eta_1$.  \\
\noindent Let us prove that
\begin{equation}\label{conv eta}
\lim\limits_{\al\nearrow 1}\eta_\al=\eta_1.
\end{equation}
\noindent Firstly, from Proposition 2 it holds that
\begin{equation}\label{c.u.}
\lim\limits_{\al\nearrow 1} f_\al(x) =f_1(x),
\end{equation}
 where the convergence is uniform over compact sets. \\
Secondly, analysing $f'_1$  we have that $f'_1(0)=1,$  $f'_1(+\infty)=0^-$, there exists a unique  
$ \eta_0\approx 0.3195$ such that $ f'_1(\eta_0)=0$ and  $f'_1(x)<0 $ for all $x>\eta_0$. In fact, $\eta_0$ is the unique positive solution to equation $ \sqrt{\pi} x e^{x^2} {\rm erfc}(x)=4x^2$. 
%\begin{equation}\label{f'_1}
% \exists\,!\, \eta_0\approx 0.3195\, \text{ such that } f'_1(\eta_0)=0, \, f'_1(x)<0 \, \forall \,x>\eta_0. 
%\end{equation}
Being $\eta_1\approx 0.6201$ it follows that $f'_1(\eta_1)<0$. Then, there exists an interval $[\eta_1-\rho, \eta_1+\rho]$, for some $\rho>0$ where $f_1$ is decreasing.\\
Finally, let $\varepsilon>0$ be ($\varepsilon<\rho$). Let $r$ be the line of equation $y=x$. Clearly  $P_1(\eta_1,\eta_1)\in r$ and we can take $P_a(a,a)$ and $P_b(b,b)$ in $r$ ($a<\eta_1<b$) such that 
\begin{equation}\label{dist}
d(P_1,P_a)<\varepsilon,\quad  d(P_1,P_b)<\varepsilon\quad \text{and}\quad  f_1 \quad \text{ is decreasing in}\, [a,b].
\end{equation}
 \noindent Being $f_1(\eta_1)=\eta_1$, it holds that $f_1(a)-a>0$ and $f_1(b)-b<0$.\\
Now let $h_0=\min\left\{f_1(a)-a, b-f_1(b)\right\}>0$. From (\ref{c.u.}) it results that there exists $\al_0 \in (0,1)$ such that 
$$
|f_\al(x)-f_1(x)|<h_0\quad \text{ for all}\, x\, \in \, [a,b], \, \text{ for all}\, \al\, \in \, (\al_0,1].
$$
Then, if $\al\in (\al_0,1]$ we have that
$$
f_\al(a)>f_1(a)-h_0>a\quad \text{and} \quad f_\al(b)<f_1(b)+h_0<b.
$$
\noindent  Applying Bolzano's Theorem ($f_\al$ is continuous in $\bbR^+$ for all $\al \in (0,1]$) it holds that the unique solution $\eta_\al $ to equation $f_\al(x)=x$ belongs to  $(a,b)$. From $(\ref{dist})$ and calling $P_\al(\eta_\al,\eta_\al)$  we get that \linebreak $|\eta_\al-\eta_1|<d(P_\al,P_1)<d(P_a,P_1)<\varepsilon $ for all $\al\,\in \, (\al_0,1]$ and (\ref{conv eta}) holds.

 %\begin{figure}[h]
%\includegraphics[width=10cm]{eta-alpha-conv-a-eta-1.eps}
%\caption{Figure 1}
%\end{figure}
 \noindent Finally, applying Propositions 1 and 2 we get that  
\begin{equation*}\begin{split}
\displaystyle\lim_{\al \nearrow 1}w_\al &= \displaystyle\lim_{\al \nearrow 1}1-\frac{1}{1-W\left(-2\eta_\al,-\frac{\al}{2},1\right)}\left[1-W\left(-\frac{x}{t^{\al/2}},-\frac{\al}{2},1\right)\right]\\
& = 1-\frac{1}{erf(\eta_1)}erf\left( \frac{x}{2\sqrt{t}} \right)
\end{split}
\end{equation*}

and
$$
\displaystyle\lim_{\al \nearrow 1}r_\al=\displaystyle\lim_{\al \nearrow 1}2\eta_\al t^{\al/2}=2\eta_1\sqrt{t}
$$
\noindent which proves that solution (\ref{u-C}) of problem (\ref{(P1)}) converges to solution (\ref{u-clasica}) of problem (\ref{St-Clasico}) as we wanted to see. 
\noindent The second part of the proof is analogous.
\endproof

\section{Conclusions}
We have considered two fractional Stefan problems involving    Riemann-Liouville and Caputo derivatives  of order $\al\in (0,1)$ such that  in the limit case ($\al=1$) both problems coincide with the same  classical Stefan problem, and the relationship between the governed equations and the Stefan conditions is  analysed.  For both problems, explicit solutions were presented and it has been proved that these solutions are different, and so, the fractional Stefan conditions are different (unless $\al=1$). Finally, the convergence when $\al \nearrow 1$ was computed obtaining for both problems the same classical solution.   

\section{Acknowledgments}
\noindent This paper has been partially sponsored by the Project PIP No. 0275 from CONICET-UA (Rosario, Argentina) and AFOSR-SOARD Grant FA 9550-14-1-0122. 

\bibliographystyle{plain}

\bibliography{Roscani_BIBLIO_GENERAL_nombres_largos2017_08}

 \end{document}